\documentclass[amscd,amssymb,11pt]{amsart}

\usepackage[all]{xy}

\newtheorem{thm}{Theorem}[section]
\newtheorem{lem}[thm]{Lemma}
\newtheorem{cor}[thm]{Corollary}
\newtheorem{prop}[thm]{Proposition}

\theoremstyle{definition}

\newtheorem{defn}[thm]{Definition}
\newtheorem{defns}[thm]{Definitions}

\newtheorem{ex}[thm]{Example}

\theoremstyle{remark}

\newtheorem{rem}[thm]{Remark}
\newtheorem{rems}[thm]{Remarks}

\numberwithin{equation}{section}

\newcommand{\thmref}[1]{Theorem~\ref{#1}}

\newcommand{\secref}[1]{\S\ref{#1}}

\newcommand{\propref}[1]{Proposition~\ref{#1}}
\newcommand{\lemref}[1]{Lemma~\ref{#1}}

\newcommand{\exref}[1]{Example~\ref{#1}}

\newcommand{\hocolim}{\operatorname*{hocolim}}

\newcommand{\colim}{\operatorname*{colim}}

\newcommand{\Map}{\operatorname{Map}}
\newcommand{\MapS}{\operatorname{Map_{\mathcal S}}}
\newcommand{\MapT}{\operatorname{Map_{\mathcal T}}}
\newcommand{\MapA}{\operatorname{Map_{\mathcal Alg}}}

\newcommand{\C}{{\mathcal  C}}
\newcommand{\Spp}{{\mathcal  S}}
\newcommand{\Sppu}{{\mathcal  S}_u}
\newcommand{\Set}{{\mathcal  Set}}
\newcommand{\T}{{\mathcal  T}}
\newcommand{\Tu}{{\mathcal  T}_u}

\newcommand{\E}{{\mathcal  E}}

\newcommand{\Ab}{{\mathcal Ab}}
\newcommand{\Alg}{{\mathcal Alg}}
\newcommand{\Algu}{{\mathcal Alg}_u}
\newcommand{\Algg}{{\mathcal Alg}^{\prime}}
\newcommand{\GT}{{\mathcal T^{\Gamma}}}

\newcommand{\Sppi}{{\text{SP}^{\infty}}}
\newcommand{\ootimes}{{\otimes}}

\newcommand{\Z}{{\mathbb  Z}}
\newcommand{\N}{{\mathbb  N}}
\newcommand{\R}{{\mathbb  R}}
\newcommand{\F}{{\mathbb  F}}
\newcommand{\Q}{{\mathbb  Q}}

\newcommand{\Sinfty}{\Sigma^{\infty}}

\newcommand{\sm}{\wedge}
\newcommand{\smb}{\overline{\wedge}}

\newcommand{\da}{\downarrow}

\newcommand{\ra}{\rightarrow}
\newcommand{\xra}{\xrightarrow}

\newcommand{\hra}{\hookrightarrow}

\begin{document}

\title[A model for a space tensored with a ring spectrum]{The McCord model for the tensor product of a space and a commutative ring spectrum}             

\author[Kuhn]{Nicholas J.~Kuhn}                                  
\address{Department of Mathematics \\ University of Virginia \\ Charlottesville, VA 22903}    
\email{njk4x@virginia.edu}
\thanks{This research was partially supported by a grant from the National Science Foundation}     
 
\date{January 28, 2002.}


\begin{abstract}   We begin this paper by noting that, in a 1969 paper in the Transactions, M.C.McCord introduced a construction that can be interpreted as a model for the categorical tensor product of a based space and a topological abelian group.  This can be adapted to Segal's very special $\Gamma$--spaces - indeed this is roughly what Segal did - and then to a more modern situation: $K \otimes R$ where $K$ is a based space and $R$ is a unital, augmented, commutative, associative $S$--algebra.  

The model comes with an easy-to-describe filtration.  If one lets $K = S^n$, and then stabilize with respect to $n$, one gets a filtered model for the Topological Andr\'e--Quillen Homology of $R$.  When $R = \Omega^{\infty} \Sinfty X$, one arrives at a filtered model for the connective cover of a spectrum $X$, constructed from its $0^{th}$ space.

Another example comes by letting $K$ be a finite complex, and $R$ the $S$--dual of a finite complex $Z$.  Dualizing again, one arrives at G.~Arone's model for the Goodwillie tower of the functor sending $Z$ to $\Sinfty \MapT(K,Z)$.

Applying cohomology with field coefficients, one gets various spectral sequences for deloopings with known $E_1$--terms.  A few nontrivial examples are given.

In an appendix, we describe the construction for unital, commutative, associative $S$--algebras not necessarily augmented.

\end{abstract}
                                                  
\maketitle

\section{Introduction} \label{introduction}

The point of this paper is to tell a story that begins with a 1969 paper of M.C.McCord \cite{mccord}, and ends with various disparate objects of current interest, e.g. Goodwillie towers, Topological Hochschild Homology, and  Topological Andr\'e--Quillen Homology.  Line by line, I think most of this story is known.  However, taken as a whole, I think the older work sheds some light on the newer.  Moreover, in this era of tremendous activity in homotopical algebra of various sorts, it seems important to remind ourselves that the genesis of many of the most useful ideas lies way back in the literature.

Conceptually, we feature the following categorical notion.  Let $\T$ be the category of based topological spaces.  If $\C$ is a category enriched over $\T$, there is the notion of the tensor product of a space $K \in \T$ with an object $X \in \C$: this is an object $K \otimes X \in \C$ satisfying
$$ \Map_{\C}(K \otimes X, Y) = \Map_{\T}(K, \Map_{\C}(X,Y))$$
for all $X,Y \in \C$.

We give an overview of the paper.

We consider various topological categories of structured objects.  Let  $\Ab$ be the category of abelian topological monoids.  Let $\GT$ be G.Segal's category of $\Gamma$--spaces \cite{segal}: functors $X: \Gamma \ra \T$, where $\Gamma$ is the category of finite based sets. Let $\Alg$ be the category of commutative, associative, augmented $S$--algebras, where $S$ is the sphere spectrum.  

These categories are closely related.  First of all, an abelian topological monoid $A$ defines in a natural way $A^{\times} \in \GT$.  In Segal's terminology, $A^{\times}$ is an example of a `special' object,  and the abelian topological groups define `very special objects'.  A very special $\Gamma$--space is roughly the same thing as an infinite loopspace, and, in this introduction, we will tempt fate and often identify these two notions. Now we note that if $X$ is either an abelian topological monoid or an infinite loopspace, then $\Sinfty X_+$ is in the category $\Alg$.

With $\C$ equal to any of these three categories, McCord's construction yields a functor\footnote{McCord uses the notation $B(A, K)$ where we use $\Sppi(K,A)$; I have borrowed my notation from \cite[Chapter 9]{ff}.}
$$ \Sppi: \T \times \C \ra \C.$$
His construction generalizes the infinite symmetric product construction studied by Dold and Thom in the 1950's: with $\N$ denoting the natural numbers, $\Sppi(K, \N) = \Sppi(K)$.

If $\C$ is either $\Ab$ or $\Alg$, there is an isomorphism 
\begin{equation} \label{tensor equation} \Sppi(K,X) = K \otimes X.
\end{equation}
See \propref{abprop}(1), and \propref{algprop}. The latter proposition (or at least its proof) seems to be new.

Equation (\ref{tensor equation}) is almost true if $\C = \GT$. S. Schwede \cite{schwede}, following Bousfield and Friedlander \cite{bf}, defines a model category structure on $\GT$, having the very special $\Gamma$--spaces as the fibrant objects, and such that $ho(\GT)$ is equivalent to the homotopy category of connective spectra.  With this structure, $\Sppi(K, \hspace{.1in} )$ preserves fibrant objects and the natural map $K \otimes X \ra \Sppi(K,X)$ is a weak equivalence.

McCord's interest stemmed from the following fundamental property: 
\begin{equation} \label{class space equation} \text{ for $X \in \Ab$, $\Sppi(S^1, X)$ is a classifying space for $X$.}
\end{equation}
Suitably interpreted, the same result is true if $X$ is  a special object of $\GT$.  For $R \in \Alg$, $\Sppi(S^1,R)$ is also of interest: $\Sppi(S^1,R) = S^1 \otimes R$ equals $THH(R;S)$, the Topological Hochschild Homology of $R$ with coefficients in the bimodule $S$.  See \propref{THHprop}; this is deduced from a similar result due to J.~McClure, R.~Schw\"{a}nzl, and R.~Vogt \cite{msv}.

The construction has two other basic properties, both discussed by McCord when $\C = \Ab$.  

Firstly, there are natural isomorphisms
\begin{equation} \Sppi(K \wedge L, X) = \Sppi(K, \Sppi(L,X)).
\end{equation}
From this and (\ref{class space equation}) one quickly deduces that $\Sppi(S^n, X)$ is an $n$--fold delooping of $X$ for $X \in \Ab$, or for very special $X \in \GT$.  For $R \in Alg$, $\Sppi(S^n,R)$ can be interpreted as `higher' Topological Hochshild Homology of $R$.

Secondly, $\Sppi(K,X)$ comes with a nice increasing filtration.  When $\C = \Ab$, one easily sees that there is an equivalence
\begin{equation} \label{filtration} F_d\Sppi(K,X)/F_{d-1}\Sppi(K,X) \simeq K^{(d)} \sm_{\Sigma_d} X^{\sm d},
\end{equation}
and little variants of this hold when $\C = \GT$ or $\Alg$.
Here $K^{(d)}$ denotes the $\Sigma_d$--space obtained from the $d$--fold smash product $K^{\sm d}$ by collapsing the the fat diagonal to a point.\footnote{This notation, which the author likes, comes from \cite{arone}.}

This much of the story will be fleshed out in sections \ref{mccord}, \ref{gamma space}, and \ref{ring spectra}, with the statement about $THH$ appearing in  \secref{taq}.

In \cite{basterra}, it is observed that $\Alg$ is Quillen equivalent to the category $\Algg$, of commutative, nonunital $S$--algebras.  In \secref{reduced section}, we discuss the corresponding filtered $\Sppi$ construction, which again agrees with the tensor product.  This reduced construction is `smaller' than what is done in $\Alg$, and the category $\E$ of finite sets and epimorphisms replaces $\Gamma$.

The quotients of the filtration, and the use of $\E$, may look vaguely familiar to readers of \cite{arone}, and this is what we explain in \secref{arone tower}.  There is a contravariant functor from $\T$ to $\Algg$ sending a based space $Z$ to the ring spectrum $R = D(Z)$, where $D$ denotes the S--dual, and the multiplication on $R$ is induced by the diagonal on $Z$.  We note that there is a natural map in $\Algg$ 
$$  K \otimes D(Z) \ra D(\MapT(K,Z))$$
and reinterprete the main convergence theorem of \cite{arone} as saying that the adjoint of this, 
$$ \Sinfty \MapT(K,Z) \ra D(K \otimes D(Z)),$$
is an equivalence if both $Z$ and $K$ are finite dimensional complexes, and the dimension of $K$ is less than the connectivity of $Z$. See \thmref{K tensor DZ}. Since $K \otimes D(Z)$ has a nice increasing filtration, the S--dual is a tower of fibrations.  This tower is visibly equivalent to the tower found by Arone, and is thus the Goodwillie tower associated to the functor sending a space $X$ to the spectrum $\Sinfty \MapT(K,Z)$.  

In \secref{taq}, we discuss the following construction.  Given $R \in \Alg$, nice properties of $\Sppi(K,R)$ as a functor of $K$ allow us to define a filtered spectrum $TAQ(R)$ by
$$ TAQ(R) = \hocolim_{n \ra \infty} \Omega^{n} \Sppi(S^n,R).$$
This is one of various equivalent definitions of the Topological Andr\'e--Quillen spectrum of $R$.  When $R = \Sinfty X_+$, with $X$ an infinite loopspace, $TAQ(R)$ is the connective delooping of $X$.  When $R = D(Z_+)$, the filtered spectrum $TAQ(R)$ is related to constructions studied by the author in \cite{k2}.

Using (\ref{filtration}), one can identify the quotients of the filtration of $TAQ(R)$:
\begin{equation} \label{TAQfiltration} F_dTAQ(R)/F_{d-1}TAQ(R) \simeq \Sigma K_d \sm_{h\Sigma_d} (R/S)^{\sm d},
\end{equation}
where $K_d$ is the $d^{th}$ partition complex which arose in the work of Arone and Mahowald on the Goodwillie tower of the identity \cite{am}.  

In \secref{examples} we note that applying ordinary cohomology with field coefficients $\mathbb F$ to the filtered spectra $\Sppi(S^n,R)$ and $TAQ(R)$ yield spectral sequences converging to $H^*(S^n \otimes R; \mathbb F)$ and $H^*(TAQ(R); \mathbb F)$, and having $E_1$ terms isomorphic to known functors of $H^*(R;\mathbb F)$: see \thmref{ss theorem}.  These spectral sequences appear to be unexplored, even in the case when $R = \Sinfty X_+$, with $X$ an infinite loop space, so that, e.g., the $TAQ(R)$ spectral sequence is calculating the cohomology of a connective spectrum from knowledge of the cohomology of its $0^{th}$--space.  As examples, we use results of ours from \cite{k2} to explain how the spectral sequence works, when $\mathbb F$ has characteristic 2, in the cases $R = \Sinfty(\Z/2_+)$, $\Sinfty(S^1_+)$,  and, most interestingly, $D(S^1_+)$.

In the Appendix, we note how our models need to be slightly tweaked when one considers commutative unital $S$--algebras not necessarily augmented.

Very influential to me in my understanding of the older work surveyed in this paper was Chapter of 9 of the unpublished book {\em The actions of the classical small categories of topology} by Bill and Ed Floyd \cite{ff}.  Writing this book was Ed's project in the late 1980's, after returning to ordinary academic life after finishing a term as provost of the University of Virginia.  I have also benefited from conversations with Mike Mandell, Bill Dwyer, and Greg Arone.

Versions of this work were presented at talks at the Johns Hopkins topology conference in the spring of 2000, and at the Union College topology and category theory conference of fall 2001.

\section{McCord's construction} \label{mccord} 

Let $K$ be a based space with basepoint $*$, and let $A$ be an abelian topological monoid.  

Imagine using (\ref{tensor equation}) to guide the construction $\Sppi(K,A)$.  As a first experiment, suppose $A = \N$.  In this case (\ref{tensor equation}) tells us that, for all $A \in \Ab$, we should have
$$ \Map_{\Ab}(\Sppi(K, \N), A) = \Map_{\T}(K, \Map_{\Ab}(\N,A)).$$
But, since $\N$ is the free abelian (topological) monoid on one generator, we can identify $\Map_{\Ab}(\N,A)$ with $A$.  Thus we are asking that $\Sppi(K,\N)$ satisfy
$$ \Map_{\Ab}(\Sppi(K, \N), A) = \Map_{\T}(K, A).$$
In other words, $\Sppi(K, \N)$ should be the free topological abelian monoid generated by $K$, a.k.a. $\Sppi(K)$.

Note that elements in $\Sppi(K)$ are words of the form $k_1^{n_1}\cdots k_d^{n_d}$, with $k_i \in K$ and $n_i \in \N$.  This suggests a definition.

\begin{defn}  Let $\Sppi(K,A)$ be the abelian topological monoid with generators $k^a$ with $k \in K$ and $a \in A$ subject to the relations: \\

\noindent (i)  $*^a = *$ for all $a \in A$, \\

\noindent (ii)  $k^0 = *$ for all $k \in K$, \\

\noindent (iii)  $k^{a_1}k^{a_2} = k^{a_1 + a_2}$ for all $a_1, a_2 \in A$. \\

\noindent Viewing $k^a$ as an element in $K \times A$, $\Sppi(K,A)$ is topologized as the evident quotient space of $\coprod_{d=0}^{\infty} (K \times A)^d$. \\
\end{defn}

Note that the abelian topological monoid satisfying only the relations of type (i) and (ii) is $\Sppi(K \sm A)$.  If $B$ is another abelian topological monoid, a monoid  map $\Sppi(K \sm A) \ra B$ corresponds to a map of topologial spaces $\phi: K \sm A \ra B$ which itself corresponds to a map $K \ra \MapT(A,B)$. This latter map takes values in $\Map_{\Ab}(A,B)$ exactly when $\phi(k,a_1+a_2) = \phi(k,a_1)\phi(k,a_2)$.  Thus we see that the quotient of $\Sppi(K\sm A)$ having type (iii) relations imposed satisfies the universal property of $K \otimes A$. (Compare with \cite[p.164]{ff}.)  We have checked the first part of the next proposition.

\begin{prop} \label{abprop} There are the following natural identifications in $\Ab$. \\

\noindent{(1)} $\Sppi(K,A) = K \otimes A$. \\

\noindent{(2)} $\Sppi(S^0,A) = A$. \\

\noindent{(3)} $\Sppi(K \sm L, A) = \Sppi(K, \Sppi(L,A))$. \\

\noindent{(4)} $\Sppi(K \vee L, A) = \Sppi(K,A) \times \Sppi(L,A)$.

\end{prop}

The last three parts of this proposition follow formally from statement (1).  For example, (3) follows by manipulating adjunctions:
\begin{equation*}
\begin{split}
\Map_{\Ab}(\Sppi(K \sm L, A), B) & = \Map_{\T}(K \sm L, \Map_{\Ab}(A,B)) \\
& = \Map_{\T}(K, \Map_{\T}( L, \Map_{\Ab}(A,B))) \\
& = \Map_{\T}(K, \Map_{\Ab}(\Sppi(L,A),B)) \\
& = \Map_{\Ab}(\Sppi(K, \Sppi(L, A)), B).
\end{split}
\end{equation*}
For statement (4), one needs to also note that the coproduct in $\Ab$ is the product. (In the next section, we will see that there are also reasonable direct proofs of (3) and (4).)

In the introduction to \cite{mccord}, McCord comments that $\Sppi( \ , A)$ ``has a tendancy to convert cofibrations \dots to quasifibrations'', and proves this in various cases \cite[Thm.8.8]{mccord}.  Note that statement (4) of the last proposition nicely illustrates his statement.

In particular, when applied to the cofibration
$ S^0 \hra I \ra S^1, $
his observations suggest that 
$\Sppi(I,A) \ra \Sppi(S^1,A)$ is a quasifibration with homotopy fiber $A = \Sppi(S^0,A)$.  One has

\begin{prop} If $A$ has a nondegenerate base point, then 
$\Sppi(S^1, A)$ is a classifying space for $A$.  
\end{prop}

Combined with statement (3) of the previous proposition, this implies

\begin{cor}  In this case, $\Sppi(S^n,A)$ is an $n$--fold classifying space of $A$.
\end{cor}

With such a mild point set hypothesis, this proposition and corollary occur as \cite[Cor.9.16]{ff}.

We end this section by noting that $\Sppi(K,A)$ is filtered by letting the $d^{th}$ filtration, which we denote $F_d\Sppi(K,A)$, be the evident quotient of $\coprod_{r=0}^{d} (K \times A)^r$.  Let $i: \Delta(K,d) \hra K^{\sm d}$ denote the inclusion of the fat diagonal into the $d$--fold smash product.  Under reasonable conditions, e.g. if $K$ is a based C.W.~ complex, the inclusion $i$ will be a $\Sigma_d$--equivariant cofibration, and we let $K^{(d)} = K^{\sm d}/\Delta(K,d)$.  Assuming the basepoint in $A$ is also nondegenerate, the inclusion $F_{d-1} \Sppi(K,A) \hra F_d\Sppi(K,A)$ will be a cofibration, and there is a homeomorphism
\begin{equation} \label{Afiltration} F_d\Sppi(K,A)/F_{d-1}\Sppi(K,A) \simeq K^{(d)} \sm_{\Sigma_d} A^{\sm d}.
\end{equation}
Statements similar to this appear in \cite[\S 6]{mccord}. \\

\begin{rem}  One should note that specializing the filtration on $\Sppi(K,A)$ to $\Sppi(K, \N) = \Sppi(K)$ does {\em not} yield the standard filtration on $\Sppi(K)$.  For example, in this paper an element $k^2 \in \Sppi(K)$ would be in filtration 1.
\end{rem}

\section{$\Gamma$--spaces and Segal's theorem} \label{gamma space}

Our first goal in this section is to rewrite our construction $\Sppi(K,A)$ in a way allowing for generalization.

We begin by defining more precisely the category $\Gamma$.

\begin{defn} Let $\Gamma$ be the category with objects the based finite sets $\mathbf  0 = \emptyset_+ $ and  ${\mathbf  n} = \{1,2,\dots, n\}_+$, $n \geq 1$, and with all based functions as morphisms.  Note that $\mathbf 0$ is both an initial and terminal object. 
\end{defn}

We remark that, unfortunately, it was the {\em opposite} of this category that was called $\Gamma$ in \cite{segal}, and the literature is strewn with inconsistent notation.

A $\Gamma$--space is then defined to be a covariant functor $X: \Gamma \ra \T$ that is `based' in the sense that it sends $\mathbf 0$ to the one point space. These are the objects of a category $\T^{\Gamma}$ having the natural transformations as morphisms.  This is a category enriched over $\T$: the set of morphisms between two $\Gamma$--spaces, $\Map_{\Gamma}(X,X^{\prime})$, has a natural topology.  Similarly, a based contravariant functor $Y: \Gamma^{op} \ra \T$ will be called a $\Gamma^{op}$--space, and these are objects in a topological category $\T^{\Gamma^{op}}$.

\begin{ex} \label{Atimes ex} If $A$ is an abelian topological monoid, there is an associated $\Gamma$--space $A^{\times}$ defined as follows.  Firstly,  $A^{\times}(\mathbf n) = A^n$.  Then, given $\alpha: \mathbf n \ra \mathbf m$, the $i^{th}$ component of $\alpha_*: A^n \ra A^m$ sends $(a_1, \dots , a_n)$ to $\prod a_j$, with the product running over $j$ such that $\alpha(j) = i$.  This product is interpreted to be the unit of $A$ if there are no such $j$.
\end{ex}

\begin{ex} Since, for any $K \in \T$, $K^n = \MapT(\mathbf n, K)$, there is the evident $\Gamma^{op}$--space $K^{\times}$ with $K^{\times}(\mathbf n) = K^n$.
\end{ex}

Note that the constructions in these last examples embed $\Ab$ into $\T^{\Gamma}$, and $\T$ into $\T^{\Gamma^{op}}$, as a full subcategories.

We now recall the coend construction.  If $X$ is a $\Gamma$--space and $Y$ is a $\Gamma^{op}$--space, we let $Y \sm_{\Gamma} X \in \T$ denote the quotient space
$$ \bigvee_n Y(\mathbf n) \sm X(\mathbf n)/(\sim),$$
where $\alpha^*(y) \sm x \sim y \sm \alpha_*(x)$ generates the equivalence relation.

It is useful to observe that, because $X(\mathbf 0) = * = Y(\mathbf 0)$, $Y \sm_{\Gamma} X \in \T = Y \times_{\Gamma} X$, where $Y \times_{\Gamma} X$ is the quotient space 
$$ \coprod_n Y(\mathbf n) \times X(\mathbf n)/(\sim),$$
where $\alpha^*(y) \sm x \sim y \sm \alpha_*(x)$ generates the equivalence relation.

By inspection, one observes 

\begin{lem} $\Sppi(K,A) = K^{\times} \sm_{\Gamma}A^{\times}$.
\end{lem}

This suggests a generalization of our construction.

\begin{defn}  Given $K \in T$ and $X \in \T^{\Gamma}$, let $\Sppi(K,X) = K^{\times} \sm_{\Gamma} X$.
\end{defn}

A couple of remarks are now in order.  

Firstly, a Yoneda's lemma type argument shows that
$$ \Sppi(\mathbf n, X) = X(\mathbf n).$$
Thus, as a functor of $K$, $\Sppi(K,X)$ extends $X$ to $\T$; more precisely, this is the left Kan extension \cite[Chap.X]{maclane}.

Secondly, we can extend our construction to 
$$ \Sppi: \T \times \T^{\Gamma} \ra \T^{\Gamma}$$
by first letting $X_n(\mathbf m) = X(\mathbf{nm})$, and then by defining $$\Sppi(K,X)(\mathbf n) = \Sppi(K, X_n).$$

It is natural to wonder if $\Sppi(K,X)$ can be then  be interpreted as a tensor product in $\T^{\Gamma}$.  Alas, this is not the case: the simple minded construction $K \sm X$ defined by $(K \sm X)(\mathbf n) = K \sm X(\mathbf n)$ is easily seen to play this role\footnote{One can then formally deduce that, given $A, B \in \Ab$, the natural map $K \sm A^{\times} \ra \Sppi(K,A^{\times})$ induces a homeomorphism $$\Map_{\Gamma}(\Sppi(K,A^{\times}), B^{\times}) = \Map_{\Gamma}(K \sm A^{\times}, B^{\times}).$$}.

Since $\Sppi(K,X)$ is not the tensor product in $\T^{\Gamma}$, formal arguments used to prove the last three statements in \propref{abprop} don't apply.  However, we still can prove suitable versions of these.

First of all, the two remarks above (along with the observation that $\mathbf 1 = S^0$) combine to show that $\Sppi(S^0, X) = X$.

Less obvious are the other two.  Statement (3) of \propref{abprop} is unchanged in our greater generality.

\begin{prop} \label{smashprop} $\Sppi(K \sm L, X) = \Sppi(K, \Sppi(L,X))$.
\end{prop}

To generalize statement (4), let $a: \Gamma \times \Gamma \ra \Gamma$ be the functor sending $(\mathbf m, \mathbf n)$ to $\mathbf {m+n}$.  Pulling back by $a$ defines $a^*: \T^{\Gamma} \ra \T^{\Gamma \times \Gamma}$.

\begin{prop} \label{wedgeprop} $\Sppi(K \vee L, X) = (K^{\times} \times L^{\times}) \times_{\Gamma \times \Gamma} a^*(X)$.
\end{prop}

\propref{abprop}(4) follows from this, once one observes that
$$ a^*(A^{\times}) = A^{\times} \times A^{\times},$$
so that
\begin{equation*}
\begin{split}
\Sppi(K \vee L, A^{\times}) & = (K^{\times} \times L^{\times}) \times_{\Gamma \times \Gamma} (A^{\times} \times A^{\times}) \\
& = (K^{\times} \times_{\Gamma} A^{\times}) \times (L^{\times} \times_{\Gamma} A^{\times}) \\
& = \Sppi(K,A^{\times}) \times \Sppi(L,A^{\times}).
\end{split}
\end{equation*}
 
Let $m: \Gamma \times \Gamma \ra \Gamma$ be the functor sending $(\mathbf m, \mathbf n)$ to $\mathbf {mn}$\footnote{More precisely, $m$ is the smash product followed by the lexicographic identification of $\mathbf m \sm \mathbf n$ with $\mathbf {mn}$.}.  We let $a_*: \T^{\Gamma^{op} \times \Gamma^{op}} \ra \T^{\Gamma^{op}}$ and $m_*: \T^{\Gamma^{op} \times \Gamma^{op}} \ra \T^{\Gamma^{op}}$ respectively denote the left adjoints to pulling back by $a$ and $m$.  We have two fundamental lemmas.

\begin{lem} \label{smashlem} $m_*(K^{\times} \times L^{\times}) = (K \sm L)^{\times}$.
\end{lem}

\begin{lem} \label{wedgelem} $a_*(K^{\times} \times L^{\times}) = (K \vee L)^{\times}$.
\end{lem}

Assuming these for the moment,  \propref{smashprop} and \propref{wedgeprop} follow.  For example, using \lemref{smashlem}, we have identifications
\begin{equation*}
\begin{split}
\Sppi(K, \Sppi(L,X)) & = K^{\times} \times_{\Gamma} (L^{\times} \times_{\Gamma} X_*) \\
& = (K^{\times} \times L^{\times}) \times_{\Gamma \times \Gamma} m^*(X) \\
& = m_*(K^{\times} \times L^{\times}) \times_{\Gamma} X \\
& = (K \sm L)^{\times} \times_{\Gamma} X \\
& = \Sppi(K \sm L, X),
\end{split}
\end{equation*}
and \propref{smashprop} follows.  The proof of \propref{wedgeprop} is similar.

\begin{proof}[Proof of \lemref{smashlem}] Given a $\Gamma^{op} \times \Gamma^{op}$--space $Y$, $m_*(Y)$ is explicitly given by 
$$m_*(Y)(\mathbf c) = \colim_{\mathbf c \da m} Y$$
where $\mathbf c \da m$ is the category with objects all triples $(\mathbf a, \mathbf b, \gamma)$ with $\gamma:\mathbf c \ra \mathbf{ab}$, and morphisms given by pairs $(\alpha: \mathbf a \ra \mathbf {a^{\prime}}, \beta: \mathbf b \ra \mathbf {b^{\prime}})$ making an appropriate diagram commute.

One such triple is $(\mathbf c, \mathbf c, \Delta)$, where $\Delta: \mathbf c \ra \mathbf {cc}$ is the diagonal, and there is a canonical morphism from this triple to any other triple $(\mathbf a, \mathbf b, \gamma)$ given by the two components of $\gamma:\mathbf c \ra \mathbf{ab}$.  It follows that the canonical map $Y(\mathbf c, \mathbf c) \ra \colim_{\mathbf c \da m} Y$ is a quotient map.

Now we specialize to $Y = K^{\times} \times L^{\times}$.  Both  $\colim_{\mathbf c \da m} K^{\times} \times L^{\times}$ and $(K \sm L)^c$ are quotients of $K^c \times L^c$, thus we just need to verify that each maps to the other, as quotient spaces of $K^c \times L^c$.

To construct a map from the former to the latter, we observe that that the smash product construction, 
$$ \MapT(A,K) \times \MapT(B,L) \ra \MapT(A \sm B, K \sm L),$$
specializes to give natural maps
$$ \sm: K^a \sm L^b \ra (K \sm L)^{ab}.$$
Thus, associated to a triple $(\mathbf a, \mathbf b, \gamma)$, there is a canonical map
$$ K^a \sm L^b \xra{\sm} (K \sm L)^{ab} \xra{\gamma^*} (K \sm L)^c,$$
and these induce the needed map $\colim_{\mathbf c \da m} K^{\times} \times L^{\times} \ra (K \sm L)^c$.

To construct a map in the other direction, we observe that $(K \sm L)^c$ is the quotient of $K^c \times L^c$ obtained by collapsing to a point the subspace 
$$\{ (k_1,\dots,k_c, l_1, \dots,l_c) \ | \ \text{ for all $i$, either $k_i = *$ or $l_i = *$} \}.$$
One checks easily that this subspace is precisely the union of the images of maps 
$$ \alpha^* \times \beta^*: K^a \times L^b \ra K^c \times L^c$$
such that the composite $\mathbf c \xra{\Delta} \mathbf {cc} \xra{\alpha \beta} \mathbf {ab}$ is the constant map $\mathbf 0$.  Thus the subspace maps to the basepoint in $\colim_{\mathbf c \da m} K^{\times} \times L^{\times}$, i.e. the projection $K^c \times L^c \ra \colim_{\mathbf c \da m} K^{\times} \times L^{\times}$ factors through $(K \sm L)^c$.
\end{proof}

\begin{proof}[Sketch proof of \lemref{wedgelem}] This follows easily from the observation that there are canonical decompositions
$$ (K \vee L)^c = \bigvee K^a \vee L^b,$$
with the wedge running over bijections $\gamma:\mathbf c \ra (\mathbf a + \mathbf b)$ which are order preserving when restricted to $\gamma^{-1}((\mathbf a + \mathbf b) - \mathbf a)$ and $\gamma^{-1}((\mathbf a + \mathbf b) - \mathbf b)$.
\end{proof}

\begin{rems} The two lemmas include the statements that wedge and smash are the left Kan extensions to $\T \times \T$ of the composites $\Gamma \times \Gamma \xra{a} \Gamma \hra \T$ and $\Gamma \times \Gamma \xra{m} \Gamma \hra \T$.

We suspect that \lemref{wedgelem} has been observed by others.   \lemref{smashlem}  seems less familiar. (Compare our proof of \propref{smashprop} to the proof of \cite[Lemma 3.7]{segal}.)  We note that $\T^{\Gamma} \times \T^{\Gamma} \xra{\sm} \T^{\Gamma \times \Gamma} \xra{m_*} \T^{\Gamma}$ is the smash product of \cite{lydakis}. \\
\end{rems}

The category of $\Gamma$--spaces admits products in the obvious way: if $X$ and $Y$ are $\Gamma$--spaces, one lets $(X \times Y)(\mathbf n) = X(\mathbf n) \times Y(\mathbf n)$.  We have

\begin{prop}  $\Sppi(K,X \times Y) = \Sppi(K,X) \times \Sppi(K,Y)$.
\end{prop}

To prove this, we first note that $X \times Y = \Delta^*(X \times Y)$, where $\Delta^*: \T^{\Gamma \times \Gamma} \ra \T^{\Gamma}$ is induced by the diagonal $\Delta: \Gamma \ra \Gamma \times \Gamma$.  Thus we have identifications
\begin{equation*}
\begin{split}
\Sppi(K, X \times Y) & = K^{\times} \times_{\Gamma} \Delta^*(X \times Y) \\
& = \Delta_* (K^{\times}) \times_{\Gamma \times \Gamma} (X \times Y) \\
& = (K^{\times} \times K^{\times}) \times_{\Gamma \times \Gamma} (X \times Y) \\
& = (K^{\times} \times_{\Gamma} X)  \times (K^{\times} \times_{\Gamma} Y) \\
& = \Sppi(K, X) \times \Sppi(K,Y),
\end{split}
\end{equation*}
where we have used the next lemma.

\begin{lem}  $\Delta_*(K^{\times}) = K^{\times} \times K^{\times}$.
\end{lem}
\begin{proof} This can be proved in various ways.  Perhaps the slickest proof is to first note that $\Delta$ is right adjoint to $a$. That $\Delta_* = a^*$ formally follows.  Finally, it is evident that $a^*(K^{\times}) = K^{\times} \times K^{\times}$.
\end{proof}

It remains, in this section, to discuss how the $\Sppi$ construction interacts with the homotopy theory of $\Gamma$--spaces.  

Define $\pi_*^s(X) = \colim_{n \ra \infty} \pi_{*+n}(\Sppi(S^n,X))$.  The colimit here arises from maps $S^1 \sm \Sppi(S^{n-1}, X) \ra \Sppi(S^n,X)$ which themselves are special cases ($K = S^1$ and $Y = \Sppi(S^{n-1},X)$) of the natural transformation
$$ K \sm Y \ra \Sppi(K , Y).$$
If we define {\em weak equivalences} to be maps $f: X \ra Y$ with $\pi_*^s(f)$ an isomorphism, then Bousfield and Friedlander \cite{bf}, following Segal \cite{segal}, showed that that the localized category $\T^{\Gamma}[weq^{-1}]$ is equivalent to the homotopy category of connective spectra.  Even more, this equivalence is induced by a Quillen equivalence between appropriate model categories.  Schwede \cite{schwede} modifies the cofibration and fibrations slightly.  All these authors work with simplicial sets rather than topological spaces, but \cite[Appendix B]{schwede} allows for some translation into our setting.

The upshot is roughly the following.  Cofibrations are maps $f:X \ra Y$ where $Y$ is obtained from $X$ by successively attaching appropriate sorts of free $\Gamma$--spaces. Fibrant objects agree with Segal's notion of a {\em very special} $\Gamma$--space, where $X$ is very special means that each map
\begin{equation} \label{special condition} X(\mathbf {a+b}) \ra X(\mathbf a) \times X(\mathbf b)
\end{equation}
is a weak equivalence of spaces, and also
\begin{equation} \text{the monoid } \pi_0(X(\mathbf 1)) \text{ is a group.}
\end{equation}

\begin{prop}  If $K$ is a C.W.~ complex, then $\Sppi(K, \hspace{.1in})$ preserves cofibrations and acyclic cofibrations.
\end{prop}

\begin{prop} If $K$ is a C.W.~ complex, and $X$ is cofibrant, then the natural map
$$ K \sm X \ra \Sppi(K,X)$$
is a weak equivalence.
\end{prop}

\begin{prop} If $K$ is a C.W.~ complex, and $X$ is cofibrant and very special, then $\Sppi(K, X)$ is again very special.
\end{prop}

\begin{thm} If $X$ is cofibrant and very special, then $\Sppi( \hspace{.1in}, X)$ takes cofibration sequences of C.W.~ complexes to a homotopy fibration sequence.  In particular, there are weak equivalences of spaces
$$X(\mathbf 1) \xra{\sim} \Omega \Sppi(S^1, X)(\mathbf 1) \xra{\sim} \Omega^2 \Sppi(S^2, X)(\mathbf 1) \xra{\sim} \dots.$$
\end{thm}

We briefly indicate why the propositions hold.  Firstly, under the cofibrancy hypotheses, $\Sppi(K,X)$ will be nicely filtered, and satisfy
\begin{equation}
F_d\Sppi(K,X)/F_{d-1}\Sppi(K,X) = K^{(d)} \sm_{\Sigma_d} (X(\mathbf d)/X_{sing}(\mathbf d))
\end{equation}
where $X_{sing}(\mathbf d)$ denotes the union of all the images of maps $X(\mathbf c) \ra X(\mathbf d)$ with $c<d$.

It follows then that then $K \sm X(\mathbf 1) \ra \Sppi(K,X)$ is a weak equivalence through a stable range, and the first two of the propositions easily can be deduced.  

For the next proposition, we note that, if (\ref{special condition}) holds, then 
$$ X_{a+b} \ra X_a \times X_b$$
is a {\em strict equivalence} of $\Gamma$--spaces, where a map is a strict equivalence if evaluating on any $\mathbf n$ yields a weak equivalence of spaces.  Then we have equivalences
\begin{equation*}
\begin{split}
\Sppi(K, X)(\mathbf {a+b}) & = \Sppi(K, X_{a+b}) \\
& \xra{\sim} \Sppi(K, X_a \times X_b) \\
& = \Sppi(K, X_a) \times \Sppi(K, X_b) \\
& = \Sppi(K, X)(\mathbf a) \times \Sppi(K,X)(\mathbf b),
\end{split}
\end{equation*}
showing that $\Sppi(K,X)$ again satisfies (\ref{special condition}).

For the theorem, see \cite[Prop.3.2]{segal} and \cite[Lemma 4.3]{bf}.  It follows that if $X$ is cofibrant and very special, then $X(\mathbf 1)$ is canonically weakly equivalent to an infinite loop space.  Furthermore, for any C.W.~ complex $K$, there are weak equivalences
\begin{equation}
F_d\Sppi(K,X)/F_{d-1}\Sppi(K,X) \simeq K^{(d)} \sm_{\Sigma_d} X(\mathbf 1)^{\sm d}.
\end{equation}

\section{Commutative ring spectra} \label{ring spectra}

We now show that the results of the previous sections extend nicely to the world of structured ring spectra.  

We work within the category $\Spp$, the category of $S$--modules studied in \cite{ekmm}.  Given $K \in \T$ and $X \in \Spp$, $\Sinfty K$, $K \sm X$, and $\Map(K,X)$ will denote the usual $S$--modules\footnote{What we are calling $\Map(K,X)$ is $F_S(\Sinfty K, X)$ in \cite{ekmm}.}.

Let $\Alg$ be the category of unital, commutative, associative, augmented $S$--algebras.  Thus an object in $\Alg$ is an $S$--module $R$, together with multiplication $\mu: R \sm R \ra R$, unit $\eta: S \ra R$ and counit $\epsilon: R \ra S$ satisfying the usual identities.  Morphisms preserve all structure.

This category is enriched over $\T$: given $R,Q \in \Alg$, the morphism space $\MapA(R,Q)$ is based with basepoint $R \xra{\epsilon} S \xra{\eta} Q$.  We also note that the coproduct in $\Alg$ of $R$ and $Q$ is $R \sm Q$.

As observed in \cite[\S1]{basterra}, results in \cite{ekmm} show that $\Alg$ has a topological model category structure in which weak equivalences are morphisms that are weak equivalences as maps of $S$--modules\footnote{In Basterra's notation, $\Alg$ is denoted $\mathcal C_{S/S}$.}.

We have two important sources of examples.

\begin{ex}  If $A \in \Ab$, then $\Sinfty A_+ \in \Alg$.  More generally, if $X$ is an $E_{\infty}$--space (e.g. an infinite loop space), then $\Sinfty X_+$ is naturally an object in $\Alg$. (See \cite[Ex.IV.1.10]{may} and \cite[\S II.4]{ekmm}.)
\end{ex}

\begin{ex}  Given a based space $Z$, let $D(Z_+)$ denote $\Map(Z_+,S)$.  This is an object in $\Alg$: the unit and the counit are respectively induced by $Z_+ \ra S^0$ and $S^0 \ra Z_+$, and the diagonal $\Delta: Z \ra Z \times Z$ induces the multiplication
$$ D(Z_+) \sm D(Z_+) \ra D(Z_+ \sm Z_+) \xra{\Delta^*} D(Z_+).$$
\end{ex}

Given $R \in Alg$, we let $R^{\sm}: \Gamma \ra \Spp$ denote the functor with $R^{\sm}(\mathbf n) = R^{\sm n}$ analogous to \exref{Atimes ex}.

\begin{defn}  Given $K \in \T$ and $R \in \Alg$, let $\Sppi(K,R) = K^{\times} \sm_{\Gamma} R^{\sm}$.
\end{defn}

We will momentarily see that $\Sppi(K,R)$ is again an object in $\Alg$.

Proofs from \secref{gamma space} extend immediately to prove the next proposition.

\begin{prop}  There are the following natural identifications. \\

\noindent{(1)} $\Sppi(S^0,R) = R$. \\

\noindent{(2)} $\Sppi(K \vee L, R) = \Sppi(K,R) \sm \Sppi(L,R)$. \\

\noindent{(3)} $\Sppi(K, R \sm Q) = \Sppi(K,R) \sm \Sppi(K,Q)$. \\
\end{prop}

A consequence of this proposition is that $\Sppi(K,R)$ takes values in $\Alg$, with multiplication given by the composite
$$ \Sppi(K,R) \sm \Sppi(K,R) = \Sppi(K,R \sm R) \xra{\Sppi(K,\mu)} \Sppi(K,R).$$
We note that this multiplication agrees with the composite
$$ \Sppi(K,R) \sm \Sppi(K,R) = \Sppi(K \vee K,R) \xra{\Sppi(\nabla,R)} \Sppi(K,R),$$
where $\nabla: K \vee K \ra K$ is the fold map.

With this structure, all the identifications in the last proposition are as objects in $\Alg$, and we also have the next proposition, whose proof follows from the arguments of the last section.

\begin{prop} \label{Rsmashprop} $\Sppi(K \sm L, R) = \Sppi(K, \Sppi(L,R))$.
\end{prop}

Now we check that $\Sppi(K,R)$ is the categorical tensor product in $\Alg$.  The following lemmas are easily verified, where we use the following notation: with $\mathcal C$ either $\T$ or $\Spp$, and $X$ and $Y$ functors from $\Gamma$ to $\mathcal C$, $\Map_{\C}^{\Gamma}(X,Y)$ denotes the space of natural transformations from $X$ to $Y$.

\begin{lem} For all $K,L \in \T$, $\Map_{\T}^{\Gamma}(K^{\times}, L^{\times}) = \MapT(K,L)$.
\end{lem}

\begin{lem} For all $R,Q \in \Alg$, $\Map_{\Spp}^{\Gamma}(R^{\sm}, Q^{\sm}) = \MapA(R,Q)$.
\end{lem}

\begin{prop} \label{algprop} For all $K \in \T$ and $R \in \Alg$, $\Sppi(K,R)$ is naturally isomorphic to $K \otimes R$.
\end{prop}

\begin{proof}  We check that $\Sppi(K,R)$ satisfies the universal property of the tensor.  Given $K \in \T$, and $R,Q \in \Alg$, we have
\begin{equation*}
\begin{split}
\MapA(\Sppi(K,R),Q) & = \Map_{\Spp}^{\Gamma}(\Sppi(K,R)^{\sm}, Q^{\sm}) \\
   & = \Map_{\Spp}^{\Gamma}(\Sppi(K,R^{\sm}), Q^{\sm}) \\
   & = \Map_{\Spp}^{\Gamma}(K^{\times} \sm_{\Gamma}m^*(R^{\sm}), Q^{\sm}) \\
& = \Map_{\T}^{\Gamma}(K^{\times}, \Map_{\Spp}^{\Gamma}(m^*(R^{\sm}), Q^{\sm})) \\
   & = \Map_{\T}^{\Gamma}(K^{\times}, \MapA(R^{\sm}, Q)) \\
   & = \Map_{\T}^{\Gamma}(K^{\times}, \MapA(R, Q)^{\times}) \\
   & = \Map_{\T}(K, \MapA(R, Q)). 
\end{split}
\end{equation*}
Here $m: \Gamma \times \Gamma \ra \Gamma$ is multiplication as in the last section.
\end{proof}

As before, $\Sppi(K,R)$ is naturally filtered.  Let $R/S$ denote the cofiber of $\eta:S \ra R$.  If $K$ is a C.W. complex, and $\eta$ is a cofibration, then 
the inclusion $F_{d-1} \Sppi(K,R) \hra F_d\Sppi(K,R)$ will be a cofibration, and there is an isomorphism of $S$--modules
\begin{equation} \label{Rfiltration} F_d\Sppi(K,R)/F_{d-1}\Sppi(K,R) \simeq K^{(d)} \sm_{\Sigma_d} (R/S)^{\sm d}.
\end{equation}

\section{The reduced model} \label{reduced section}

It is sometimes useful to replace $\Alg$ by a slightly different category.  Let $\Algg$ be the category of nonunital, commutative, associative $S$--algebras (the category denoted $\mathcal N_S$ in \cite{basterra}).  Basterra observes that the functor $S \vee : \Algg \ra \Alg$, that wedges a unit $S$ onto a nonunital algebra, has as right adjoint the augmentation ideal functor $J: \Alg \ra \Algg$,
defined by letting $J(R)$ be the fiber of $R \xra{\epsilon} S$.  She then notes that, with a natural topological model category on $\Algg$, these adjoint functors form a Quillen pair, and thus induce adjoint equivalences on the associated homotopy categories.

\begin{ex}  If $Z$ is a based space, $J(D(Z_+))= D(Z)$.  The multiplication on $D(Z)$ is induced by the reduced diagonal $\Delta: Z \ra Z \sm Z$.
\end{ex}

Our $\Sppi(K, \hspace{.1in})$ construction has a `reduced' analogue in $\Algg$.

\begin{defn}  Let $\E$ be the category with objects $\mathbf n$, for $n \geq 1$, and with morphisms from $\mathbf n$ to $\mathbf m$ equal to all epimorphisms from $\{1, \dots, n\}$ to $\{1, \dots, m\}$.
\end{defn} 

As observed in \cite{arone} (see also \cite{ahearnkuhn}), a based space $K$ defines a functor $K^{\sm}: \E^{op} \ra \T$ with $K^{\sm}(\mathbf n) = K^{\sm n}$.  Also, $J \in \Algg$ defines $J^{\sm}: \E \ra \Spp$ in the obvious way.

\begin{defn} Given $K \in T$ and $J \in \Algg$, let $\Sppi(K,J) = K^{\sm} \sm_{\E} J^{\sm}$.
\end{defn}

The analogues of all the properties of $\Sppi(K,R)$ proved in the last section hold in our setting, with virtually identical proofs.  In particular, $\Sppi(K,J)$ is again an object in $\Algg$, and it agrees with the categorical tensor product $K \otimes J$. 

From the above comments, one can formally deduce the following isomorphism in $\Alg$.

\begin{prop} $\Sppi(K,J) \vee S = \Sppi(K,J \vee S)$.
\end{prop}

Though we won't show this here, this proposition can also be given a direct proof, and there are analogues in other contexts.  Readers may wish to compare this result with observations in \cite{pirashvili}.

As usual, $\Sppi(K, J)$ is filtered: if $\E_d$ denote the full subcategory of $\E$ with objects $\mathbf n$ for $n \leq d$, then we let $F_d \Sppi(K,J) = K^{\sm} \sm_{\E_d} J^{\sm}$.  If $K$ is a C.W. complex,  then 
the inclusion $F_{d-1} \Sppi(K,J) \hra F_d\Sppi(K,J)$ will be a cofibration, and there is an isomorphism of $S$--modules
\begin{equation} \label{Jfiltration} F_d\Sppi(K,J)/F_{d-1}\Sppi(K,J) \simeq K^{(d)} \sm_{\Sigma_d} J^{\sm d}.
\end{equation}

We note that the isomorphism of the last proposition is filtration preserving.

\section{Reinterpretation of Arone's tower for $\Sinfty \MapT(K,X)$} \label{arone tower}

In this section, we let $K$ be a finite C.W.~ complex.  

In \cite{arone}, G.~ Arone described a model for the Goodwillie tower of the functor sending a based space $Z$ to the $S$--module $\Sinfty \MapT(K,Z)$.  Here we show that, if $Z$ is a finite complex, his tower arises as the S--dual of the filtered object $K \ootimes D(Z)$ of the last section.  

We recall Arone's construction and some of its properties \cite{arone}.  For more detail, see also \cite{ahearnkuhn}.

\begin{defns} Let $Z$ be a based space. 

\noindent{(i)} Let $P^K(Z) = \Map_{\Spp}^{\E}(K^{\sm}, Z^{\sm})$, the spectrum of natural transformations from $\Sinfty K^{\sm}$ to $\Sinfty Z^{\sm}$. 

\noindent{(ii)} Let $P_d^K(Z) = \Map_{\Spp}^{\E_d}(K^{\sm}, Z^{\sm})$. 

\noindent{(iii)}  Let $\Phi(K,Z): \Sinfty \MapT(K,Z) \ra P^K(Z)$ be the natural transformation that sends $f:K \ra Z$ to the natural transformation with $n^{th}$ component equal to $\Sinfty f^{\sm n}: \Sinfty K^{\sm n} \ra \Sinfty Z^{\sm n}$. \\
\end{defns} 

The spectrum $P^K(Z)$ is the inverse limit of the tower of fibrations
 $$\dots \ra P^K_{d+1}(Z) \ra P^K_d(Z) \ra P^K_{d-1}(Z) \ra \dots ,$$  
and the $d^{th}$ fiber is isomorphic to
$$ \Map_{\Spp}^{\Sigma_d}(K^{(d)}, Z^{\sm d}).$$
Because $K^{(d)}$ is finite, and the $\Sigma_d$ action on this is free away from the basepoint, this fiber is naturally homotopy equivalent to the homotopy orbit spectrum
$$ (D(K^{(d)}) \sm Z^{\sm d})_{h \Sigma_d}. $$

From this last description one sees that the tower has the form of a Goodwillie tower, and also that the connectivity of the fibers goes up if the connectivity of $Z$ is greater than the dimension of $K$.  Arone then proves that this {\em is} the Goodwillie tower of $\Sinfty \MapT(K,Z)$ by proving

\begin{thm} \cite{arone} If the connectivity of $Z$ is greater than the dimension of $K$, then $\Phi(K,Z)$ is a homotopy equivalence.
\end{thm}

Now we connect these constructions to $K \ootimes D(Z)$.

\begin{defns} Let $Z$ be a based space. 

\noindent{(i)} Let $\Tilde{\Theta}(K,Z):K \ootimes D(Z) \ra D(\MapT(K,Z))$
be the map in $\Algg$ adjoint to the composite
$$ K \xra{\text{eval}} \MapT(\MapT(K,Z),Z) \xra{D} \Map_{\Algg}(D(Z), D(\MapT(K,Z)).$$ 

\noindent{(ii)}  Let $\Theta(K,Z): \Sinfty \MapT(K,Z) \ra D(K \ootimes D(Z))$ be the $S$--module map adjoint to $\Tilde{\Theta}(K,Z)$.

\noindent{(iii)}  Let $\alpha(K,Z): P^K(Z) \ra D(K \ootimes D(Z))$ be the map of  cofiltered $S$--modules defined as follows.  Let $i: Z^{\sm} \ra D(D(Z)^{\sm})$ be the natural transformation adjoint to $D(Z)^{\sm} \ra D(Z^{\sm})$.  Now let
$\alpha(K,Z)$ be the composite induced by $i$:
$$ \Map_{\Spp}^{\E}(K^{\sm}, Z^{\sm}) \xra{i} \Map_{\Spp}^{\E}(K^{\sm}, D(D(Z)^{\sm})) = \MapS(K^{\sm}\sm_{\E} D(Z)^{\sm}, S).$$
\end{defns}

A check of the definitions verifies the next lemma.

\begin{lem}  There is a commutative diagram
\begin{equation*}
\xymatrix{ 
& \Sinfty \MapT(K,Z) \ar[dl]_{\Phi(K,Z)} \ar[dr]^{\Theta(K,Z)} & \\
P^K(Z) \ar[rr]^-{\alpha(K,Z)} &  & D(K \ootimes D(Z))
}
\end{equation*}
\end{lem}

\begin{lem}  $\alpha(K,Z)$ is a homotopy equivalence if $Z$ is a finite complex.
\end{lem}
\begin{proof}  If $Z$ is finite, then $i: Z^{\sm n} \ra D(D(Z)^{\sm n})$ is a an equivalence for all $n$.  Now the lemma follows by observing that $K^{\sm}$ is a cofibrant $\E^{op}$--space, or more simply, note that the fibers of the towers will be equivalences, as $K^{(d)}$ is a free $\Sigma_d$--complex for all $d$. (This has been noted before; see e.g. \cite{ahearnkuhn, mccarthy}.)
\end{proof}

Summarizing, we conclude

\begin{thm} \label{K tensor DZ} If both $K$ and $Z$ are finite complexes, and the dimension of $K$ is less than the connectivity of $Z$, then
$$\Theta(K,Z): \Sinfty \MapT(K,Z) \ra D(K \ootimes D(Z))$$
is a weak equivalence, and thus the algebra map
$$\Tilde{\Theta}(K,Z):K \ootimes D(Z) \ra D(\MapT(K,Z))$$
can be identified as the map from a spectrum to its double dual.
\end{thm}

We end this section by noting how the homological version of this discussion would go.  

Let $\mathbb F$ be a field, $H \mathbb F$ the associated commutative $S$--algebra, and $\Alg_{\F}$ the category of commutative, nonunital $H \F$ algebras.  Let $K \otimes_{\F}J \in \Alg_{\F}$ denote the tensor product of a based space $K$ and an $J \in \Alg_{\F}$. As before, one learns that 
$$ K \otimes_{\F} J = K^{\times} \sm_{\E} J^{\sm},$$
where smash products are taken over $H \F$.

Now let $D_{\F}(Z) = \Map(Z, H\F)$, the $H \F$--module whose homotopy groups are the cohomology groups of $Z$ with $\F$--coefficients.  In this case, the natural map $i: H\F \sm Z^{\sm n} \ra D_{\F}(D_{\F}(Z)^{\sm n})$
is an equivalence for any space $Z$ with $H_*(Z; \F)$ of finite type.  Reasoning as before, from Arone's theorem one deduces

\begin{thm}  If $K$ is a finite complex of dimension less than the connectivity of $Z$, and $H_*(Z; \F)$ is of finite type, then the natural map in $\Alg_{\F}$,
$$ \Theta: K \otimes_{\F} D_{\F}(Z) \rightarrow D_{\F}(\MapT(K,Z)),$$
is an equivalence.
\end{thm}

\begin{rem}  It seems likely that this theorem can be deduced from older convergence results for the Anderson spectral sequence \cite{anderson}, and then one can run our arguments backwords, and {\em deduce} Arone's theorem.  The novelty would then be to identify the filtration as Arone did.
\end{rem}

\section{Topological Hochschild homology and Topological Andre--Quillen homology} \label{taq}

Let $THH(R;M)$ denote the Topological Hochschild homology spectrum associated to a $S$--algebra $R$ and an $R$--bimodule $M$ (see e.g. \cite[Chap.9]{ekmm}).  If $R$ is commutative and augmented, then $\epsilon: R \ra S$ makes $S$ into an $R$--bimodule.  We have

\begin{prop} \label{THHprop} $S^1 \otimes R = THH(R;S)$.
\end{prop}
\begin{proof}  This is a variant of a theorem of J.~McClure, R.~Schwanzl, and R.~Vogt \cite{msv}.  They show that if $R$ is a commutative $S$--algebra, then $THH(R;R)$ is the tensor product of $R$ with $S^1$ with the tensor product in the category of commutative $S$--algebras.  In the appendix, we note that if $R$ is also augmented, then this would agree with $S^1_+ \otimes R \in \Alg$.  Thus
$$ THH(R;R) = S^1_+ \otimes R.$$
Applying $\otimes R$ to the pushout square in $\T$
\begin{equation*}
\xymatrix{
S^0 \ar[d] \ar[r] & {*} \ar[d]  \\
S^1_+ \ar[r] &  S^1,}
\end{equation*}
yields a pushout square in $\Alg$

\begin{equation*}
\xymatrix{
R \ar[d] \ar[r] &
S \ar[d]  \\
S^1_+ \otimes R \ar[r] &
S^1 \otimes R,
}
\end{equation*}
and we conclude that 
$$S^1 \otimes R = (S^1_+ \otimes R) \sm_R S = THH(R,R) \sm_R S = THH(R;S).$$
\end{proof}

Given $K \in \T$ and $R \in \Alg$, there is a natural map
$$ K \sm R \ra K \otimes R, $$
and thus 
$$ K \sm (L \otimes R) \ra K \otimes (L \otimes R) = (K \sm L) \otimes R.$$
This map is easily seen to be filtration preserving.

Specializing to the case when $K = S^1$, and $L= S^n$ yield filtration preserving maps
$$ \Sigma (S^n \otimes R) \ra S^{n+1} \otimes R,$$
or, equivalently, 
$$ S^n \otimes R \ra \Omega (S^{n+1} \otimes R).$$

\begin{defn} Let $\displaystyle TAQ(R) = \hocolim_{n \ra \infty} \Omega^n S^n \otimes R.$
\end{defn}

M.~Mandell has shown the author that this definition agrees with other definitions of Topological Andr\'e Quillen Homology in the literature, e.g. \cite{basterra}.  In particular $TAQ(R)$ is homotopy equivalent to the cofiber of $J(R) \sm J(R) \ra J(R)$\footnote{Strictly speaking, R should be replaced by by a fibrant object in $\Alg$ and then $J(R)$ replaced by a cofibrant object in $\Algg$.}.

As the next example makes clear, $TAQ(R)$ can be viewed as an `infinite delooping' of $R$.

\begin{ex} If $X$ is a connective $S$--module, $TAQ(\Sinfty (\Omega^{\infty} X)_+) \simeq X.$  To see this, just recall that $S^n \otimes $ \  yields the $(n-1)$--connected $n$--fold delooping of an infinite loopspace.
\end{ex}

Note that $TAQ(R)$ is filtered with 
$$ F_d TAQ(R)/ F_{d+1} TAQ(R)  \simeq \hocolim_{n \ra \infty} \Sigma^{-n} S^{n(d)} \sm_{\Sigma_d} (R/S)^{\sm d}.$$

As in \cite{am}, let $K_d$ be the unreduced suspension of the classifying space of the poset of nontrivial partitions of a set with $d$ elements.

\begin{lem} \cite{ad}  There is a $\Sigma_d$--equivariant map 
$$ \hocolim_{n \ra \infty} \Sigma^{-n} S^{n(d)} \ra \Sigma K_d $$
that is a nonequivariant equivalence.
\end{lem}

The original short proof of this, due to Arone and Mahowald, appears in \cite[Appendix]{k2}.

\begin{cor} There is a homotopy equivalence
$$ F_d TAQ(R)/ F_{d+1} TAQ(R)  \simeq (\Sigma K_d \sm (R/S)^{\sm d})_{h\Sigma_d}.$$
\end{cor}

\section{Spectral sequences and examples} \label{examples}

Applying homology or cohomology with $\F$--coefficients to our filtered models for $S^n \otimes R$ and $TAQ(R)$ yields highly structured convergent spectral sequences with $E_1$ terms equal to  known functors of $H_*(R;\F)$.  To see why this is true, we note that there is an explicit equivariant duality map \cite{ahearnkuhn}
$$ F(\R^n,d)_+ \sm S^{n(d)} \ra S^{nd},$$
where $F(\R^n, d)$ is the usual configuration space of $d$ distinct points in $\R^n$.  Thus the homology calculations of \cite{clm} apply.

To be more precise, let $\{E_r^{*,*}(S^n \otimes R;\F)\}$ and $\{E_r^{*,*}(TAQ(R);\F)\}$ respectively denote the spectral sequences for computing $H^*(S^n \otimes R\F)$ and $H^*(TAQ(R);\F)$.  Let $\Tilde H^*(R; \F)$ denote the reduced cohomology of $R$, i.e. $H^*(J(R); \F)$.

\begin{thm} \label{ss theorem} For $R \in \Alg$ with $H_*(R; \F)$ of finite type, there are natural isomorphisms as follows.
\begin{enumerate}
\item If $\F$ has characteristic 0, then 
$$ E_1^{*,*}(S^n \otimes R;\F) = S^*(\Sigma^{1-n}L(\Sigma^{-1}\Tilde H^*(R;\F))), $$
and 
$$ E_1^{*,*}(TAQ(R);\F) = \Sigma L(\Sigma^{-1}\Tilde H^*(R;\F))).$$
\item
If $\F$ has characteristic p, then 
$$ E_1^{*,*}(S^n \otimes R;\F) = S^*(\mathcal R_n( \Sigma^{1-n}L_r(\Sigma^{-1}\Tilde H^*(R;\F)))), $$
and 
$$ E_1^{*,*}(TAQ(R);\F) = \mathcal R( \Sigma L_r(\Sigma^{-1}\Tilde H^*(R;\F)))).$$
\end{enumerate}
\end{thm}
In this theorem, $\Sigma^d V$ denotes the $d$--fold shift of a graded vector space $V$, $L$ is the free Lie algebra functor, $L_r$ is the free restricted Lie algebra functor, $S^*$ is the free commutative algebra functor\footnote{One has sign conventions of the usual sort.}, and $\mathcal R$ and $\mathcal R_n$ are appropriate free Dyer--Lashof operation functors.

The author plans to write more about this elsewhere.

We end with three examples.  All are nontrivial, and most of what I say follows immediately from work done in \cite{k2}.  More detail about the last example will appear in \cite{k5}.

\begin{ex} If $R = \Sinfty S^1_+$, then $TAQ(R) \simeq \Sigma H \Z$.  At least when localized at 2, the filtration of $TAQ(R)$ will correspond to the symmetric product of spheres filtration of $H \Z$\footnote{This presumably holds integrally: \cite[Thm.1.14]{ad} says that the filtration quotients are correct.}.  In particular, the $d^{th}$ associated graded spectrum is contractible unless $d$ is a power of $2$,  
$$ F_{2^k} = \Sigma SP^{2^k}(S), \text{ and } F_{2^k}/F_{2^{k-1}} = \Sigma^{k+1} L(k).$$
Here $L(k)$ is as in \cite{mp}.
The boundary maps of the filtration yields the complex of spectra
$$ \cdots \ra \Sigma L(2) \ra \Sigma L(1) \ra \Sigma L(0) \ra H\Z$$
occurring in the Whitehead conjecture \cite{k1}.  The sequence is exact in homotopy, but zero in mod 2 homology: indeed, the spectral sequence for computing $H_*(\Sigma H \Z; \F_2) = \Sigma A/ASq^1$ collapses at $E_1$.
\end{ex}

\begin{ex} If $R = \Sinfty \Z/2_+$, then $TAQ(R) \simeq H \Z/2$.  The $d^{th}$ associated graded spectrum is contractible unless $d$ is a power of $2$,  
$$ F_{2^k} = SP^{2^k}_{\Delta}(S^0), \text{ and } F_{2^k}/F_{2^{k-1}} = \Sigma^k M(k).$$
Here we recall \cite{mp} that $SP^{2^k}_{\Delta}(S^0)$ is defined to be the cofiber of the `diagonal' $\Delta: SP^{2^{k-1}}(S^0)\ra SP^{2^k}(S^0)$, and $M(k) = L(k) \vee L(K-1)$.
As in the previous example, the boundary maps of the filtration yields the complex of spectra
$$ \cdots \ra M(2) \ra M(1) \ra M(0) \ra H\Z/2$$
occurring in the mod 2 Whitehead conjecture \cite{k}.  The sequence is exact in homotopy, but the spectral sequence for computing $H_*(H \Z/2; \F_2) = A$ collapses at $E_1$.
\end{ex}

\begin{ex} If $R = D(S^1_+)$, then $TAQ(R) \simeq \Sigma^{-1}H\Q$.  There are various ways to see this; in \cite{k5}, we will prove that $S^2 \otimes R \simeq \Sigma H\Q \vee S^0$.  Localized at 2, the $d^{th}$ associated graded spectrum is contractible unless $d$ is a power of $2$, and 
$$ F_{2^k}/F_{2^{k-1}} = \Sigma^{-1}SP^{2^k}_{\Delta}(S^0).$$
Thus  $H^*(F_{2^k}/F_{2^{k-1}}; \F_2) = \Sigma^{-1} A/L_{k+1}$, where $L_{k}$ is the span of all admissible sequences in the Steenrod algebra of length at least $k$.  The boundary maps of the filtration yields a complex of spectra
$$ \dots \ra \Sigma^{-2}SP^{4}_{\Delta}(S^0) \ra \Sigma^{-1}SP^{2}_{\Delta}(S^0)\ra S^0$$
that is exact in cohomology: each map sends the bottom class of the cyclic $A$--module to $Sq^1$ applied to the bottom class of the next module.  
\end{ex}

\appendix

\section{Augmented versus nonaugmented ring spectra}

Let $\Algu$ be the category of unital commutative $S$--algebras, but not necessarily augmented.  Thus we have forgetful maps
$$ \Alg \ra \Algu \ra \Algg.$$

$\Algu$ is enriched over $\Tu$, the category of {\em unbased} topological spaces, so one can look for a convenient model for $K \otimes R$ with $K \in \Tu$ and $R \in \Algu$.  In this appendix we describe such a model, and compare this to the construction in \secref{ring spectra}.

\subsection{$K \otimes R$ for unital commutative algebras}

Let $\Sppu$ be the category of $S$--modules under $S$, so an object is an $S$--module map $\eta: S \ra X$.

\begin{defn} Given $K \in \Tu$ and $X \in \Sppu$, let $K \smb X \in \Sppu$ be the pushout:
\begin{equation*}
\xymatrix{
K_+ \sm S \ar[d] \ar[r] &
S \sm S = S \ar[d]  \\
K_+ \sm X \ar[r] &
K \smb X.
}
\end{equation*}
\end{defn}

It is easy to check 
\begin{lem}  There is an adjunction
$$ \Map_{\Sppu}(K \smb X, Y) = \Map_{\Tu}(K, \Map_{\Sppu}(X,Y)).$$
\end{lem}

Let $\Set$ be the category of finite sets.  Given $K \in \Tu$, there is an apparent functor $K^{\times}: \Set^{op} \ra \Tu$.  Now note that, as it was in $\Alg$, $\sm$ is the coproduct in $\Algu$.  Then $R \in \Algu$ defines $R^{\sm}: \Set \ra \Sppu$.

\begin{lem} For all $K,L \in \Tu$, $\Map_{\T}^{\Set}(K^{\times}, L^{\times}) = \Map_{\Tu}(K,L)$.
\end{lem}

\begin{lem} For all $R,Q \in \Algu$, $\Map_{\Sppu}^{\Set}(R^{\sm}, Q^{\sm}) = \Map_{\Algu}(R,Q)$.
\end{lem}

\begin{defn} Given $K \in \Tu$ and $R \in \Algu$, let 
$$\Sppi(K,R) = K^{\times} \smb_{\Set} R^{\sm}.$$
\end{defn}

As in \secref{ring spectra}, the lemmas combine to prove the analogue of \propref{algprop}.

\begin{prop} \label{algpropu} For all $K \in \Tu$ and $R \in \Algu$, $\Sppi(K,R)$ is again in $\Algu$ and is naturally isomorphic to the categorical tensor product $K \otimes R$.
\end{prop}

$\Sppi(K,R)$ is filtered in the usual way, and one gets an isomorphism of $S$--modules
\begin{equation}
F_d\Sppi(K,R)/F_{d-1}\Sppi(K,R) \simeq (K_+)^{(d)} \sm_{\Sigma_d} (R/S)^{\sm d}.
\end{equation}
We note that $(K_+)^{(d)}$ is just $S$ if $d=0$ and $K^{\times d}/(\text{fat diagonal})$ if $d>0$.

\subsection{The unbased versus the based construction}

In this subsection, we denote the tensor in $\Alg$ by $\otimes$ and the tensor in $\Algu$ by $\otimes_u$.

Given $R \in \Algu$, the product $R \times S$ will be in $\Alg$, with augmentation given by projection, and unit $S \xra{\Delta} S \times S \ra R \times S$. This construction is right adjoint to the forgetful functor:

\begin{lem} Given $Q \in \Alg$ and $R \in \Algu$, there is an adjunction
$$ \Map_{\Algu}(Q,R) = \Map_{\Alg}(Q, R \times S).$$
\end{lem}

Note that, if $Q \in \Alg$ and $R \in \Algu$, then $\Map_{\Algu}(Q,R)$ is based with basepoint $Q \ra S \ra R$.

\begin{prop}  Given $K \in \T$, $Q \in \Alg$, and $R \in \Algu$, there is an adjunction isomorphism
$$ \Map_{\Algu}(K \otimes Q, R) = \Map_{\T}(K, \Map_{\Algu}(Q,R)).$$
\end{prop}

\begin{proof}  We have natural isomorphisms
\begin{equation*}
\begin{split}
\Map_{\Algu}(K \otimes Q, R) & = \Map_{\Alg}(K \otimes Q, R \times S) \\
   & = \Map_{\T}(K, \Map_{\Alg}(Q,R \times S)) \\
   & = \Map_{\T}(K, \Map_{\Algu}(Q,R)). 
\end{split}
\end{equation*}
\end{proof}

\begin{cor} If $Q \in \Alg$ and $L \in \Tu$, then $L_+ \otimes Q = L \otimes_u Q$.
\end{cor}

\begin{proof}  Let $K = L_+$ in the proposition, and note that
\begin{equation*}
\begin{split}
\Map_{\T}(L_+, \Map_{\Algu}(Q,R)) 
   & = \Map_{\Tu}(L, \Map_{\Algu}(Q,R)) \\
   & = \Map_{\Algu}(L \otimes_u Q, R). 
\end{split}
\end{equation*}
\end{proof}

\end{document}